\documentclass[letterpaper, 10 pt, conference]{ieeeconf}  

\IEEEoverridecommandlockouts                              
\newcommand{\eq}{\begin{equation}}
\newcommand{\eeq}{\end{equation}}
\newcommand{\eqn}{\begin{eqnarray}}
\newcommand{\eeqn}{\end{eqnarray}}

\newcommand{\bsea}{\begin{subeqnarray}}
\newcommand{\esea}{\end{subeqnarray}}
\newcommand{\nn}{\nonumber}

\newtheorem{remark}{Remark}[section]
\newtheorem{teor}{Theorem}[section]
\newtheorem{corr}{Corollary}[section]
\newtheorem{propo}{Proposition}[section]
\newtheorem{lemm}{Lemma}[section]
\newtheorem{exam}{Example}
\newtheorem{obss}{Observation}
\newtheorem{probl}{Problem}
\newtheorem{defn}{Definition}[section]
\newcommand{\teo}{\begin{teor}}
\newcommand{\eteo}{\end{teor}}
\newcommand{\cor}{\begin{corr}}
\newcommand{\ecor}{\end{corr}}
\newcommand{\pro}{\begin{propo}}
\newcommand{\epro}{\end{propo}}
\newcommand{\lemma}{\begin{lemm}}
\newcommand{\elemma}{\end{lemm}}
\newcommand{\ex}{\begin{exam}}
\newcommand{\eex}{\end{exam}}
\newcommand{\pb}{\begin{probl}}
\newcommand{\epb}{\end{probl}}
\newcommand{\df}{\begin{defn}}
\newcommand{\edf}{\end{defn}}
\newcommand{\aprop}{\begin{apropo}}
\newcommand{\eaprop}{\end{apropo}}
\newcommand{\alem}{\begin{alemm}}
\newcommand{\ealem}{\end{alemm}}
\newcommand{\rema}{\begin{remark}}
\newcommand{\erema}{\end{remark}}
\newcommand{\oss}{\begin{obss}}
\newcommand{\eoss}{\end{obss}}

\usepackage{times} 
\usepackage{amsmath} 
\usepackage{amssymb}  
\usepackage{amsmath}
\usepackage{amssymb}
\usepackage{makeidx}
\usepackage{multicol}
\usepackage[bottom]{footmisc}
\usepackage{graphicx}
\usepackage{color}
\usepackage{epstopdf}
\usepackage{url}
\usepackage{algorithm}
\usepackage{algcompatible}
\usepackage{calc}

\newcommand{\Pc}{ \mathcal{P}}

\newcommand{\Sc}{ \mathcal{S}}

\newcommand{\Es}{ \mathbb{E}}

\newcommand{\Rs}{ \mathbb{R}}

\newcommand{\Zs}{ \mathbb{Z}}

\newcommand{\yv}{\mathrm{y}}
\newcommand{\ev}{\mathrm{e}}

\newcommand{\pp}{\mathbf{p}}

\title{\LARGE \bf
A Bayesian Approach to Sparse plus Low rank Network Identification}

\author{Mattia Zorzi and Alessandro Chiuso
\thanks{This work has been partially supported by the FIRB project ``Learning
meets time'' (RBFR12M3AC) funded by MIUR.}
\thanks{M. Zorzi is with the Dipartimento di Ingegneria dell'Informazione, Universit\`a degli studi di
Padova, via Gradenigo 6/B, 35131 Padova, Italy
        {\tt\small zorzimat@dei.unipd.it}}%
\thanks{A. Chiuso is with the Dipartimento di Ingegneria dell'Informazione, Universit\`a degli studi di
Padova, via Gradenigo 6/B, 35131 Padova, Italy
        {\tt\small chiuso@dei.unipd.it}}%
}

\begin{document}

\maketitle
\thispagestyle{empty}
\pagestyle{empty}

\begin{abstract}
We consider the problem of modeling multivariate stochastic processes with parsimonious dynamical models which can be represented   with a sparse dynamic network with few latent nodes. This structure translates into a sparse plus low rank model. In this paper,
we propose a Bayesian approach to identify such models.
\end{abstract}

\section{INTRODUCTION}

This paper deals with the network identification problem of a multivariate stochastic process of high dimension.
Important applications can be found in many fields, for instance in econometrics, \cite{STOCK_MARKET,Chandrasekaran_latentvariable,MULTIRES_CHOI_WILLSKY_2010,LATENTG}, social network analysis \cite{SOCIAL_NETWORK}, system biology, \cite{SYSTEM_BIOLOGY} and so on.

Let $y$ be a stochastic process whose components (or variables) $y_k$ are manifest, i.e. 
 can be directly measured. 
We define as network of this manifest process $y$ a directed graph wherein nodes denote the variables and  edges encode conditional Granger causality relations among these variables, \cite{GRANGER_CAUSALITY}. More specifically, there is an edge from node $i$ (i.e.  variable $y_i$) to node $j$ (i.e. variable $y_j$) if the past of $y_i$ is needed to predict $y_j$ conditioned on the past of all the other $y_k$, $k\neq i$. Sparse graphs (i.e. with few edges) represent a concise and interpretable way to describe the relations among the variables of the manifest process, \cite{CHIUSO_PILLONETTO_SPARSE_2012}. However, this modeling assumption does not exploit nor encode the fact that the manifest process may often be thought  as  driven by a low dimensional latent process, i.e. a stochastic process whose components cannot be directly measured.

In this paper we formulate a new network identification problem which takes into account the presence of this low dimensional latent process and the variables of the manifest process Granger causes each other mostly through the latent variables. The corresponding network has a two layer structure: one layer denotes the manifest variables and the other one the latent variables. The presence of few latent variables should drastically reduce the edges in the manifest layer, therefore increasing the degree of conciseness and robustness of the model. Finally, it turns out that a model described by this network has a sparse plus low rank (S+L) structure. More precisely, the low rank part depends on the number of latent variables, whereas the sparse one on the number of conditional Granger causality relations among the manifest variables.

When it comes to developing and identification algorithm, which maps measured data into an estimated dynamical model, it is fair to say that  the prediction error method (PEM) is a consolidate paradigm in system identification \cite{LJUNG_SYS_ID_1999,SaddleOM_STOICA_1988}. In the traditional setting, candidate models are described in fixed parametric model structures, e.g. ARMAX, whose complexity is determined using cross validation or information based criteria.   Regularization has been recently  introduced in the PEM framework, see
\cite{PILLONETTO_2011_PREDICTION_ERROR,PILLONETTO_DENICOLAO2010,EST_TF_REVISITED_2012,KERNEL_METHODS_2014}, as an alternative approach to control complexity of the estimated models. This latter class of methods
 start with  a large enough (in principle infinite dimensional) model class; the inverse (ill-posed) problem of determining a specific model from a finite set of measured data can be made into a well posed problem  using a penalty term,  whose duty is to select models with specific features. In the Bayesian view,  this is equivalent to the introduction of an {\em a priori} probability (i.e. prior) on the model to estimate. For instance, the prior should account the fact that the model is  Bounded Input Bounded Output (BIBO) stable,  \cite{PILLONETTO_DENICOLAO2010}.

The identification algorithm we propose belongs to this latter class of methods, where  the predictor impulse response is estimated in the framework of Gaussian regression. The predictor impulse response is   modeled as a zero mean Gaussian random vector. Its covariance matrix, referred to as kernel matrix, encodes the {\em a priori information}. In this case, the {\em a priori} information is that the predictor impulse responses are BIBO stable and that the model has a S+L structure. As we will see, such {\em a priori} information can be encoded in the kernel matrix using the maximum entropy principle. Moreover, this kernel matrix is characterized by the decay rate of the predictor impulse responses, by the number of conditional Granger causality relations among the manifest variables, and by the number of latent variables. These features are tuned by the so called hyperparameters  vector. It is estimated by minimizing the negative log-likelihood of the measured data. Beside the fact that this problem is nonconvex, the joint estimation of the hyperparameters tuning the sparse and low rank part is not trivial, because these two parts may be nonidentifiable from the measured data. We propose an algorithm to estimate these hyperparameters which imposes and ``hyper-regularizer'' on the low rank hyperparameter to handle partially this non-uniqueness. Once  the kernel matrix is fixed, an unique estimate of the S+L model is guaranteed through regularization.

We warn the reader that the present paper only reports some preliminary result regarding the Bayesian estimation of S+L models.  In particular, all the proofs and most of the technical assumptions needed therein are omitted and will be published afterwards.

The outline of the paper follows. 
In Section \ref{section_SL}, we introduce the S+L model for multivariate stochastic processes. In Section \ref{section_pb_formulation}, we apply  the Gaussian approach to S+L system identification. In Section \ref{section_kernel}, we derive the kernel matrix by the maximum entropy principle. In Section \ref{section_hyperparameters}, we present the algorithm to estimate the hyperparameters  of the kernel matrix. In Section  \ref{section_simulation}, we provide some numerical examples to show the effectiveness of our method. Finally, the conclusions are in Section
\ref{section_conclusions}.

\subsection*{Notation} Throughout the paper, we will use the following notation. 
$\Sc_m^+$ denotes the cone of the positive definite symmetric matrices, and $\overline{\Sc}_m^+$ its closure. Given $v\in\Rs^m$ and $G\in\Rs^{n\times m}$, $v_i$ denotes the $i$-th entry of $v$ and $[G]_{ij}$ denotes the entry of $G$ in position $(i,j)$. $\|G\|_Q$ denotes the weighted Frobenius norm of $G$ with weight matrix $Q\in\Sc_m^+$. 
Given a transfer matrix $L(z)$ of dimension $m\times m$, with some abuse of terminology, we say that $L(z)$ has rank equal to $n$, with $n\leq m$, if it admits the decomposition $L(z)=FH(z)$ where $F\in\Rs^{m\times n}$ and $H(z)$ is a $n\times m$ transfer matrix. 
Given a stochastic process $y=\{y(t)\}_{t\in\Zs}$, with some abuse of notation, $y(t)$
will both denote a random vector and its sample value.
Finally, \eq \yv^-(t):=\left[
                                 \begin{array}{ccc}
                                  y(t-1)^\top  & y(t-2)^\top & \ldots  \\
                                 \end{array}
                               \right]^\top\eeq
denotes the past data vector of $y$ at time $t$. In similar way, $\yv_i^-(t)$ denotes the past data vector of $y_i$.

\section{S+L Models} \label{section_SL}
Consider two zero mean stationary Gaussian stochastic processes $y$ and $x$ of dimension $m$ and $n$, respectively.
Let $y$ be manifest, i.e. it can be measured, and $x$ latent, i.e. it cannot be measured. We assume that $y(t) \in \Rs^m$ and $x (t)\in \Rs^n$ are described by the model
\eqn  \label{yx}y(t)&=& Fx(t)+S(z)y(t)+v(t)\nn\\
 x(t)&=&H(z)y(t)+w(t)\eeqn
where $H(z)=\sum_{k=1}^\infty H_k z^{-k}$  is a BIBO stable transfer matrix  of dimension $n\times m$, $F\in\Rs^{m\times n}$, $S(z)=\sum_{k=1}^\infty S_k z^{-k}$ is a  BIBO stable transfer matrix of dimension $m\times m$, $v$ and $w$  are, respectively, $m$ and $n$ dimensional  white Gaussian noise (WGN) with zero mean,  and covariance matrix $\Sigma_v$ and $\Sigma_w$. Moreover, $v$ and $w$ are independent.

\begin{figure}[htbp]
\begin{center}
\includegraphics[width=\columnwidth]{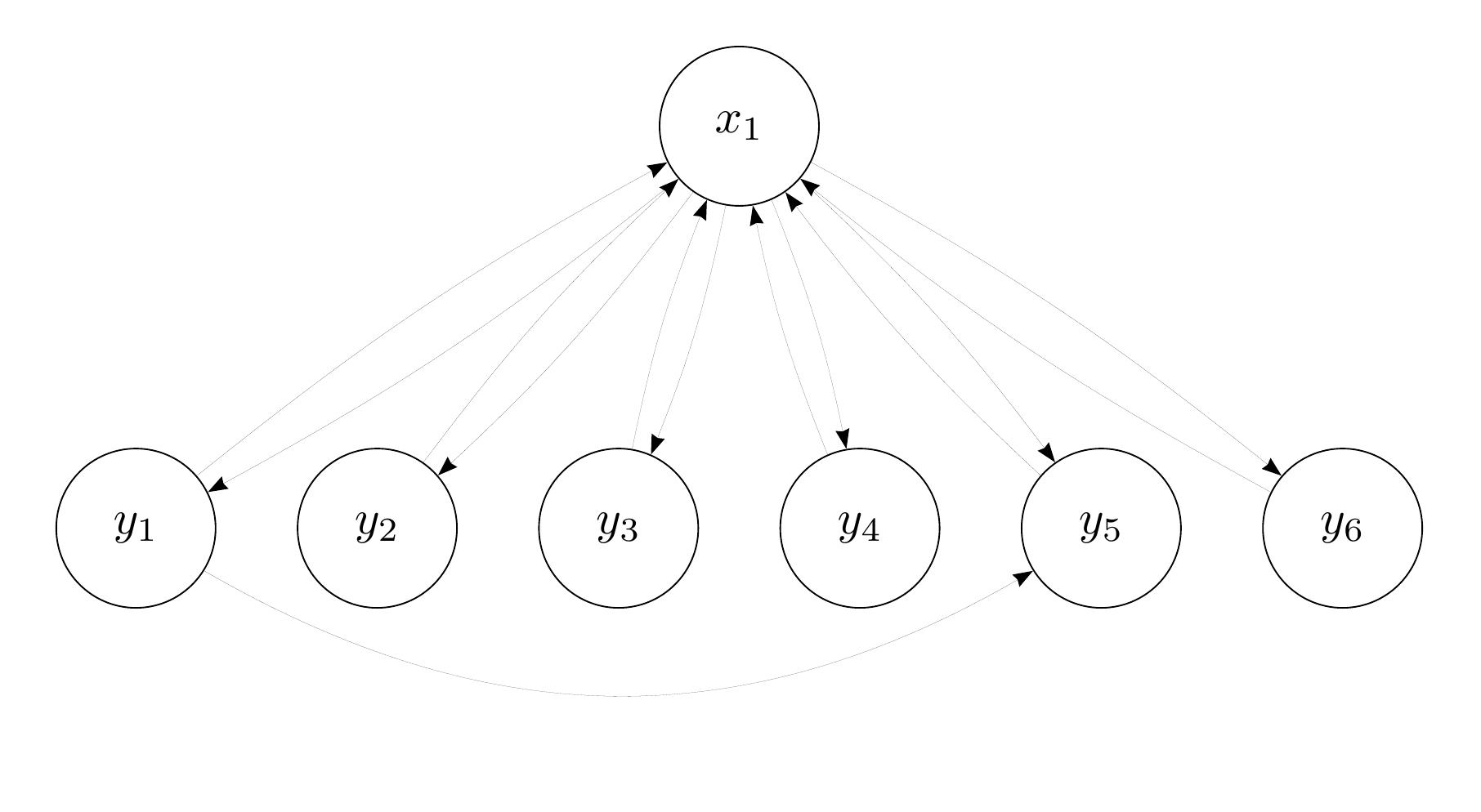}
\end{center}
\caption{Example of a sparse plus low rank model with $m=6$ and $n=1$.}\label{graph_ex}
\end{figure}

It is possible to describe the structure of model  (\ref{yx}) using a directed graph (i.e.
network) as in Figure \ref{graph_ex}, \cite{LAURITZEN_1996}. Each node of this network corresponds to components of  $y$
 or $x$; edges encode conditional Granger causality relations.  More precisely,  there is a direct link from  node  $y_j$ ($x_j$) to  node  $y_i$ ($x_i$) if and only if
 $\yv_j^-(t)$ ($x_j(t)$) is needed to predict $y_i(t)$ ($x_i(t)$), 
 conditionally on all the past information available at  time $t$, i.e. $\yv_k^-(t)$ with $k\neq j$ and $x(t)$ ($\yv^-(t)$ and $x_k(t)$ with $k\neq j$).
 In this case, we shall say that $y_j$ ($x_j$) conditionally Granger causes $y_i$ ($x_i$), \cite{GRANGER_CAUSALITY}. In Figure \ref{graph_ex} we provide an example with $m=6$ and $n=1$. In particular,  $x_1$ conditionally Granger causes all the components in $y$ and conversely, and $y_1$ conditionally Granger causes $y_5$.

Our main modeling assumption in  (\ref{yx}) is that  manifest variables Granger cause each other mostly through  few latent variables. Therefore, we have $n\ll m$, i.e. the number of latent variables in $x$ is small as compared to the number of manifest variables in $y$; similarly  $S(z)$ is
sparse, i.e. many of its entries are null transfer functions, so that   few manifest variables conditionally Granger causes each other. From (\ref{yx}), we obtain the sparse plus low rank (S+L) model for $y$:
\eq \label{S+L model} y(t)=S(z)y(t)+L(z)y(t)+e(t)\eeq
where $S(z)$ is a sparse transfer matrix by assumption, $L(z):=FH(z)$ is a low rank transfer matrix because  $F$   and $H^\top(z)$ are tall matrices, and $e(t):=v(t)+w(t)$ is WGN with covariance matrix $\Sigma=\Sigma_v+F \Sigma_w F^\top$.

Let $\hat y(t|t-1)$ be the minimum variance one-step ahead predictor of $y(t)$ based on the observations $\yv^-(t)$, and $\hat x(t|t-1)$ be the minimum variance estimator of $x(t)$ based on $\yv^-(t)$. From (\ref{yx}), we have \eqn 
 \label{yx_hat} \hat y(t|t-1)&=&F\hat x(t|t-1)+S(z)y(t)\nn\\
 \hat x(t|t-1)&=& H(z)y(t)  .\eeqn
that is,  the predictable part of $y$ is a function of few estimated latent variables and of the pasts of few manifest variables of $y$. 
Eq. \eqref{yx_hat} can be compactly written as 
 \eq 
\begin{array}{rcl} \hat y(t|t-1)&=&S(z)y(t)+L(z)y(t).
\end{array} \eeq

It is worth noting that in the case that $L(z)=0$, i.e. there is no need of latent variables to characterize the predictor of $y$, we obtain the sparse model presented in \cite{CHIUSO_PILLONETTO_SPARSE_2012}. In the case $S(z)=0$, i.e. the predictor of $y$ is completely characterized by the estimators of the latent variables, we obtain a quasi-static factor model where the noise process is white, see for instance \cite{deistler2007}. We would like to stress that  the decomposition of a transfer matrix into sparse plus low rank may not be unique. As noticed in \cite{Chandrasekaran_latentvariable}, this degeneracy may occur when $L(z)$ is sparse or $S(z)$ has few null entries. Although this nonidentifiability issue is important, the aim of this paper is to find one S+L decomposition (see Section \ref{section_pb_formulation}) which is not necessarily unique.

 \section{Gaussian Regression Approach to System Identification}\label{section_pb_formulation}
Consider model (\ref{yx}) and  assume
that the measured data $y(1) \ldots y(N)$ are extracted from a realization of $y$. The latent process $x$ cannot be measured nor its dimension $n$ is known. In this Section, we address the problem of estimating    $S(z)$ and $L(z)$ from the given data. We draw inspiration from the  Gaussian regression approach proposed in \cite{PILLONETTO_2011_PREDICTION_ERROR}. According this method, $y$ is generated by model
\eq \label{OEmodel}y(t)=G(z)y(t)+e(t)\eeq
where $G(z)$ is a BIBO stable $m\times m$ transfer matrix and 
$e$ is WGN with covariance matrix $\Sigma$. 
Note that, if we set $G(z)=S(z)+L(z)$ then (\ref{OEmodel}) is equivalent to (\ref{S+L model}).  On the other hand, with (\ref{OEmodel})
we loose the S+L structure we are interest in.

Since $G(z)$ is BIBO stable, and thus the impulse response coefficients decay to zero as a function of the lag index, it is possible to use the approximation
\eq G(z)= \sum_{k=1}^T G_k z^{-k}\eeq
where $T$ is sufficiently large. The parameters $G_k$, $k=1,.,T$ of the truncated transfer matrix  are stacked in the vector $\theta\in\Rs^{m^2 T}$ which is defined as follows
\eqn \label{def_theta} \theta &=& \left[
                   \begin{array}{ccc|c}
                    (g^{[11]})^\top & \ldots  & (g^{[1m]})^\top & \ldots \\
                   \end{array}
                 \right.\nn\\
                 && \hspace{0.2cm }\left.
                   \begin{array}{c|ccc}
                    \ldots & (g^{[m1]})^\top & \ldots &(g^{[mm]})^\top \\
                   \end{array}
                 \right]^\top
                 \eeqn
  where   \eq  \label{def_g_ij} g^{[ij]}= \left[
                  \begin{array}{cccc}
                    [G_1]_{ij} & [G_2]_{ij} & \ldots & [G_T]_{ij}\\
                  \end{array}
                \right]^\top \eeq denotes the impulse response coefficients of the transfer function in position $(i,j)$ of  the truncated approximation of $G(z)$.
Then, we stack the measured data in the vector $\yv$ as follows
  \eqn \label{y_vettore}\yv&=&\left[
          \begin{array}{ccc|c}
            y_1(T+1)^\top &  \ldots & y_1(N)^\top  & \ldots  \\
          \end{array}
        \right.\nn\\
        && \hspace{0.2cm}\left.
          \begin{array}{c|ccc}
            \ldots  & y_m(T+1) ^\top&  \ldots & y_m(N)^\top \\
          \end{array}
        \right]^\top.\eeqn In similar way, we define
  \eqn \label{e_vettore}\ev&=&\left[
          \begin{array}{ccc|c}
            e_1(T+1)^\top &  \ldots & e_1(N)^\top  & \ldots \\
          \end{array}
        \right.\nn\\ && \hspace{0.2cm} \left.
          \begin{array}{c|ccc}
            \ldots  & e_m(T+1)^\top &  \ldots & e_m(N)^\top \\
          \end{array}
        \right]^\top .\eeqn
From (\ref{OEmodel}) the vector of the measured data can be expressed in the linear regression form
\eq \yv=\Phi \theta+\ev \eeq where $\Phi\in\Rs^{mN\times m^2 T}$ is the regression matrix. Note that, 
$\Phi\theta$ is the one-step ahead predictor of $\yv$.

According to the Gaussian regression framework, $\theta$ is modeled as a zero mean Gaussian random vector with covariance matrix, or kernel matrix, denoted by $K\in\Sc_{m^2 T}^+$.

Let $\hat \theta $ be the posterior mean  of $\theta$ given $\yv$. In \cite{PILLONETTO_2011_PREDICTION_ERROR} it has been proved that, under some technical assumption, $\hat \theta$ can also be written as
solution to the Tikhonov regularization problem
\eq \label{Tikho_pb_theta}\hat \theta=\underset{\substack{\theta}}{\arg\min}\|\yv-\Phi \theta\|^2_{\Sigma^{-1}\otimes I_N}+\| \theta\|^2_{K^{-1}}.\eeq
Moreover, it is not difficult to see that \eq \label{Bayes_theta} \hat \theta= K\Phi^\top (\Phi K\Phi^\top+\Sigma\otimes I_N)^{-1}\yv.\eeq
\rema Using the theory of reproducing kernel Hilbert spaces, \cite{ARONSZAJN1950}, it is possible to show that the above results still hold for $T\rightarrow \infty$, \cite{PILLONETTO_2011_PREDICTION_ERROR}.\erema

\rema Although we assumed $K\in\Sc_{m^2 T}^+$, the Bayes estimator (\ref{Bayes_theta}) also holds for $K$ singular. In that case, Problem (\ref{Tikho_pb_theta}) is well defined provided that $\theta$ belongs
to the range of $K$.\erema

The optimal solution $\hat \theta$ highly depends on the choice of $K$. A typical assumption is that the transfer functions in $G(z)$ are independent, that is
$g^{[ij]}$ are independent vectors.  We shall also assume that $g^{[ij]}$, with $i,j=1\ldots m$, are identically distributed, so that
$K=I_{m^2}\otimes \tilde K$
where $\tilde K\in\Sc_{T}^+$ is the covariance matrix of $g^{[ij]}$, with $i,j=1\ldots m$.
In this paper, $\tilde K$ is chosen as a filtered version of the tuned/correlated (TC) kernel, see 
  \cite{PILLONETTO_2011_PREDICTION_ERROR} and \cite{EST_TF_REVISITED_2012} for more details. 
It is important to note this kernel is able to capture high frequency oscillations, which are typical of the predictor impulse responses for low pass processes, and enforces BIBO stability as $T\rightarrow \infty$ on the posterior mean of the 
predictor impulse responses.

\subsection{Gaussian Regression Approach to S+L Identification} \label{sec_SL_Gaussian_Regression}
The idea is to model $S(z)$ and $L(z)$ through Gaussian process, similarly to what has been done above for $G(z)$. 
In particular, since $S(z)$ and $L(z)$ are BIBO stable, we can consider their truncated approximations with $T$ is sufficiently large.
Then, it is not difficult to see that 
\eq \yv=\Phi(\theta_l+\theta_s)+\ev \eeq
where $\yv$ and $\ev$ have been defined in (\ref{y_vettore})  and (\ref{e_vettore}). $\theta_s,\theta_l \in\Rs^{m^2 T}$ contains the parameters of the truncated approximations of $S(z)$ and $L(z)$, respectively. Then, we model $\theta_s$ and $\theta_l$ as zero mean random vectors with covariance matrix $K_S$ and $K_L$, respectively. Moreover, we shall assume that $\theta_s$ and $\theta_l$
are Gaussian and independent. As we will see in Section \ref{section_kernel}, these assumptions are suggested by the maximum entropy principle.

\pro Let $\hat \theta_s$ and $\hat \theta_l$ be, respectively, the posterior mean of $\theta_s$ and $\theta_l$ given $\yv$. Then,  under some technical assumption,
$\hat \theta_s$ and $\hat \theta_l$ are solution to the Tikhonov regularization problem  
\eq \label{Tickho_SL} \underset{\substack{\theta_s,\theta_l}}{\arg\min}\|\yv-\Phi (\theta_s+\theta_l)\|^2_{\Sigma^{-1}\otimes I_N}+\| \theta_s\|^2_{K_S^{-1}}+\| \theta_l\|^2_{K_L^{-1}}.\eeq Moreover, we have 
\eq \hat \theta_s= K_S\Phi^\top c,\;\; \hat \theta_l=K_L\Phi^\top c\eeq
where \eq c=(\Phi(K_S+K_L) \Phi^\top+\Sigma\otimes I_N)^{-1}\yv.\eeq \epro

In what follows, $\hat \theta_s$ and $\hat \theta_l$ will be referred to as posterior mean of $S(z)$ and $L(z)$, respectively.
In the next Section, we shall show how $K_S$ and $K_L$ can be chosen 
so as  to enforce BIBO stability on both the posterior mean of $S(z)$ and $L(z)$, sparsity on the posterior mean of $S(z)$ as well as  low rank of the posterior mean of $L(z)$.

 \section{Maximum Entropy Kernel Matrix} \label{section_kernel}
In this Section we characterize the prior probability density of $\theta_s$ and $\theta_l$
by using the maximum entropy principle. Such principle states that
among all the prior probability densities satisfying certain desired constraints, the optimal one should maximize the differential entropy.

Our starting assumptions are that $\theta_s$ and $\theta_l$ are absolutely continuous zero mean  random vectors. 
Let $\pp(\theta_s,\theta_l)$ denote the joint probability density of $\theta_s$ and $\theta_l$. Let $\Es$ denote the integration over $\Rs^{2m^2 T}$ with respect to the probability measure $\pp$. Moreover, $\Pc$ denotes the space of probability densities which are Lebesgue integrable. The differential entropy of $\pp\in\Pc$ is, \cite{COVER_THOMAS},
 \eq \mathbf{H}(\pp)=-\Es[ \log(\pp(\theta_s,\theta_l))].\eeq

Next, we characterize the constraint on $\theta_s$ enforcing BIBO stability and sparsity on the posterior mean of $S(z)$.
The transfer function in position $(i,j)$ of $S(z)$ is the null transfer function if and only if $s^{[ij]}$ is the null vector.
We consider the constraint \eq \label{constraint_p_s}\Es[\| s^{[ij]}\|^2_{\tilde K^{-1}}]\leq p_{ij},\eeq
where $ p_{ij}\geq 0$.
If $p_{ij}=0$, then $s^{[ij]}$ is zero in mean square and
 so also is its posterior mean. Moreover, simple algebraic manipulations show that the weighted second moment bound in  (\ref{constraint_p_s})
 implies a bound on the variance of $k$-th element of $s^{[ij]}$ which decays as the $k$-th element in the main diagonal of $\tilde K$.
Therefore, condition (\ref{constraint_p_s}) enforces BIBO stability and sparsity on the posterior mean of $S(z)$.

Regarding the low rank constraint on $\theta_l$, let $A_l\in \Rs^{m\times m T}$ be the random matrix such that \eq A_l=\left[\begin{array}{cccc}L_1 & L_2 & \ldots & L_T\end{array}\right].\eeq
Consider the constraint 
\eq \Es[A_l (\tilde K^{-1} \otimes I_m)A_l^\top]\leq Q.\eeq
If $Q\in\overline{\Sc}_m^+$ has $m-n$ singular values equal to zero, then the posterior mean of $A_lA_l^\top$ has rank less than or equal to $n$. Therefore, the latter admits the decomposition
 \eq A_l=\left[\begin{array}{cccc}F H_1 & FH_2 & \ldots & FH_T\end{array}\right],\eeq
where $F \in \Rs^{m\times n}$ and $H_k \in \Rs^{n\times m}$, $k=1\ldots T$, as in Section \ref{section_SL}. Equivalently, the posterior mean of $L(z)$  admits the decomposition $L(z)=FH(z)$.
Similarly to the sparse part, the weight matrix $\tilde K^{-1}\otimes I_m$ enforces BIBO stability on the posterior mean of $L(z)$.

Consider the following maximum entropy problem
 \eqn  \label{ME_problem}&\underset{\pp\in\Pc}{\max}& \mathbf{H}(\pp)\nn\\
&\hbox{s.t.} & \Es[\|s^{[ij]}\|^2_{\tilde K^{-1}}]\leq p_{ij}\;\; i,j=1\ldots m\nn\\
& & \Es[A_l (\tilde K^{-1}\otimes I_m)A_l^\top]\leq Q\eeqn
where $p_{ij}>0$ $i,j=1\ldots m$, and $Q\in {\Sc}_m^+$.
\teo\label{teo_ME} The optimal solution
to (\ref{ME_problem}) is such that $\theta_s$ and $\theta_l$ are independent, Gaussian with zero mean and covariance matrix
\eq  K_S= \Gamma \otimes \tilde K,\;\;
  K_L= \Lambda \otimes I_m \otimes \tilde K\eeq
where
 \eq \Gamma=\mathrm{diag}(\gamma_{11}\ldots \gamma_{m^2} )\in\Sc_{m^2}^+\eeq and $\Lambda \in\Sc_m^+ $.\eteo

The matrices $\Gamma$ and $\Lambda$ are the hyperparameters  of the kernel matrices $K_S$ and $K_L$, respectively. 
Clearly, we are interested in the limiting case where $p_{ij}=0$ for some $(i,j)$ and $Q$ low-rank. 
Let $\mathbf{P}=\{(i,j) \hbox{ s.t } p_{ij}=0\}$ and $\mathbf{Q}=\{v\in\Rs^{m} \hbox{ s.t } Qv=0\}$.
Then, it can be shown that the maximum entropy solution can be extended by continuity to this limiting case where $\gamma_{(i-1)m+j}=0$ if and only if 
$(i,j)\in\mathbf{P}$ and $\Lambda v=0$ if and only if $v\in\mathbf{Q}$. Thus, $\Gamma$ tunes sparsity on the posterior mean of $S(z)$ and $\Lambda$ tunes the rank on the posterior mean of $L(z)$.  Finally, is worth noting that the hyperparameters  tuning the decay rate of the posterior mean of the predictor impulse responses are encoded in $\tilde K$, \cite{PILLONETTO_DENICOLAO2010}.

\section{Estimation of the Hyperparameters}\label{section_hyperparameters}
In order  to compute $\hat \theta_s$ and $\hat \theta_s$ we need to estimate the hyperparameters in $\tilde K$ and the matrices $\Gamma$ and $\Lambda$.
The hyperparameters describing  $\tilde K$ are estimated in a preliminary step by minimizing the negative log-likelihood of model (\ref{OEmodel}), see \cite{PILLONETTO_2011_PREDICTION_ERROR}.
$\Gamma$ and $\Lambda$ are obtained minimizing
the negative log-marginal likelihood $\ell$ of $\yv$.
Under  some technical assumption,  we have, \cite{PILLONETTO_2011_PREDICTION_ERROR}, 
\eq \label{loglik} \ell(\yv,\Gamma,\Lambda)=\frac{1}{2}\log \det V+\frac{1}{2}\yv^\top V^{-1} \yv+\mathrm{const. term}\eeq
where \eq \label{V_for_SL}V= \Phi (K_S+K_L)\Phi^\top+\Sigma \otimes I_{N}.\eeq
Since (\ref{loglik}) is nonconvex in $V$, only
local minima can be computed. Beside that, the joint minimization of $\Gamma$ and $\Lambda$ is not trivial because the sparse and low rank part may be nonidentifiable from the measured data. For this reason, we constrain the structure of $\Lambda$ as follows:
\eq \Lambda =\alpha (I-UU^{\top})+U\mathrm{diag} (\beta_1\ldots \beta_r) U^{\top}\eeq
where $U\in\Rs^{m\times r} $ and its columns are the first $r$ singular vectors of an estimate $\hat A_l \hat A_l^\top$ of $A_lA_l^\top$.
In this way, the constraints in $\Lambda$ are decoupled along the ``most reliable'' $r$ singular vectors of $\hat A_l \hat A_l^\top$ and their orthogonal complement. This is equivalent to fix $r$ latent variables (from the estimate $\hat A_l \hat A_l^\top$). 
Regarding the hyperparameter $\Gamma$, in \cite{ARAVKIM_CONVEX_NONCONVEX_2014} 
it has been shown that the minimization of (\ref{loglik}) leads sparsity in the main diagonal of $\Gamma$. Therefore, we minimize (\ref{loglik}) with respect to 
$\tilde \xi=\{\gamma_1\ldots \gamma_{m^2},\alpha,\beta_1\ldots \beta_r\}$ while $r$ and $U$ are fixed. The complete procedure to estimate $r$, $U$ and $\tilde \xi$ is described in Algorithm \ref{algo:hyper}. $r^{(k)}$, $U^{(k)}$, $\hat A_l^{(k)}$ and $\tilde \xi ^{(k)}$ denote, respectively, $r$, $U$, $\hat A_l$ and $\tilde \xi$ at the $k$-th iteration. Finally, to minimize efficiently (\ref{loglik}) with respect to $\tilde \xi$ we used the
scaled gradient projection algorithm developed in \cite{BONETTINI_2014}.

\begin{algorithm}
\caption{Computation of $r$, $U$ and $\tilde \xi$}
\label{algo:hyper}
\begin{algorithmic}[1] \small
\STATE $k=0$
\STATE $r^{(0)}\leftarrow 0$
\STATE $U_{OPT}^{(0)} \leftarrow $ empty matrix 
\STATE $\tilde \xi^{(0)} \leftarrow \underset{\tilde \xi}{\mathrm{argmin}}\,\ell(\yv^+,\tilde \xi,U_{OPT}^{(0)})$ 
\STATE $\tilde \xi_{OPT}^{(0)} \leftarrow \tilde \xi^{(0)}$
\REPEAT 
\STATE $k \leftarrow k+1$
\STATE $r^{(k)}\leftarrow r^{(k-1)}+1$
\IF{$r^{(k)}=1$}
\STATE $\hat A_l^{(k)}\leftarrow [\, \hat G_1\;\;  \hat G_2\;\; \ldots \,] $ where  $ \hat G_1, \hat G_2,\ldots$ are the
\STATEx \hspace{0.5cm}   coefficients of $\hat G(z)$ estimated from (\ref{OEmodel}) with 
\STATEx \hspace{0.5cm}  $K=I_{p^2}\otimes \tilde K$ 
\ELSE
\STATE $\hat A_l^{(k)}\leftarrow [\,\hat  L_1^{(k)}\;\; \hat  L_2^{(k)}\;\; \ldots \,] $ where  $ \hat L_1^{(k)}, \hat L_2^{(k)},\ldots$ are 
\STATEx \hspace{0.5cm} the coefficients of $\hat L^{(k)}(z)$ estimated from (\ref{Tickho_SL}) with 
\STATEx \hspace{0.5cm}  $K_S$ and $K_L$ having hyperparameters given by 
\STATEx \hspace{0.5cm} $U_{OPT}^{(r^{(k)}-1)}$  and $\tilde \xi_{OPT}^{(r^{(k)}-1)}$  
\ENDIF
\STATE $U^{(k)}\leftarrow$ first $r^{(k)}$ singular vectors of  $\hat  A_l^{(k)} \left. \hat A_l^{(k)}\right.^\top$
\STATE $\tilde \xi^{(k)} \leftarrow \underset{\tilde \xi}{\mathrm{argmin}}\,\ell(\yv^+,\tilde \xi,U^{(k)})$ 
\REPEAT
\STATE $\tilde \xi_{OPT}^{(r^{(k)})} \leftarrow \tilde \xi^{(k)}$
\STATE $U_{OPT}^{(r^{(k)})} \leftarrow U^{(k)}$
\STATE $k \leftarrow k+1$
\STATE $r^{(k)}\leftarrow r^{(k-1)}$ 
\STATE $\hat A_l^{(k)}\leftarrow [\, \hat L_1^{(k)}\;\; \hat  L_2^{(k)}\;\; \ldots \,] $ where  $ \hat L_1^{(k)}, \hat L_2^{(k)},\ldots$ are 
\STATEx \hspace{0.5cm} the coefficients of $\hat L^{(k)}(z)$ estimated from (\ref{Tickho_SL}) 
\STATEx \hspace{0.5cm}  with $K_S$ and $K_L$ having hyperparameters given by 
\STATEx \hspace{0.5cm}  $U_{OPT}^{(r^{(k)}-1)}$ and $\tilde \xi_{OPT}^{(r^{(k)}-1)}$  
\STATE $U^{(k)}\leftarrow$ first $r^{(k)}$ singular vectors of  $ \hat A_l^{(k)} \left. \hat A_l^{(k)}\right.^\top$
\STATE $\tilde \xi^{(k)} \leftarrow \underset{\tilde \xi}{\mathrm{argmin}}\,\ell(\yv^+,\tilde \xi,U^{(k)})$ 
\UNTIL{ $\ell(\yv^+,\tilde \xi^{(k)},U^{(k)})<\ell(\yv^+,\tilde \xi^{(k-1)},U^{(k-1)})$ }
\UNTIL{ }
\STATEx \hspace{0.2cm} $\ell(\yv^+,\tilde \xi_{OPT}^{(r^{(k)})},U_{OPT}^{(r^{(k)})})<\ell(\yv^+, \tilde \xi_{OPT}^{(r^{(k)}-1)},U_{OPT}^{(r^{(k)}-1)})$
\STATE $r  \leftarrow r^{(k)}-1$
\STATE $U\leftarrow U_{OPT}^{(r^{(k)}-1)}$
\STATE $\tilde \xi\leftarrow  \tilde \xi^{(r^{(k)}-1)}_{OPT}$
\end{algorithmic}
\end{algorithm}

\section{Simulation results}\label{section_simulation}
We consider three Monte Carlo studies of $50$ runs where at any run a model with $m=6$
manifest variables is randomly generated.
For each run in the Monte Carlo experiments an identification data set and a test set,
both of size $500$, are generated. The noise covariance matrix $\Sigma$ is always estimated via a preliminary step using a low-bias ARX-model, see \cite{GOODWIN_1992}.

In the first experiment, the models have McMillan degree equal to $20$,
and are perturbed versions of (\ref{S+L model})
with $l=1$ latent variable and four non null transfer functions in $S(z)$.

The second experiment is identical to the first one, with the exception that
the latent variables are $l=2$.

In the third experiment, the models have McMillan degree equal to $30$, but without a special structure.

\begin{figure}[htbp]
\begin{center}
\includegraphics[width=\columnwidth]{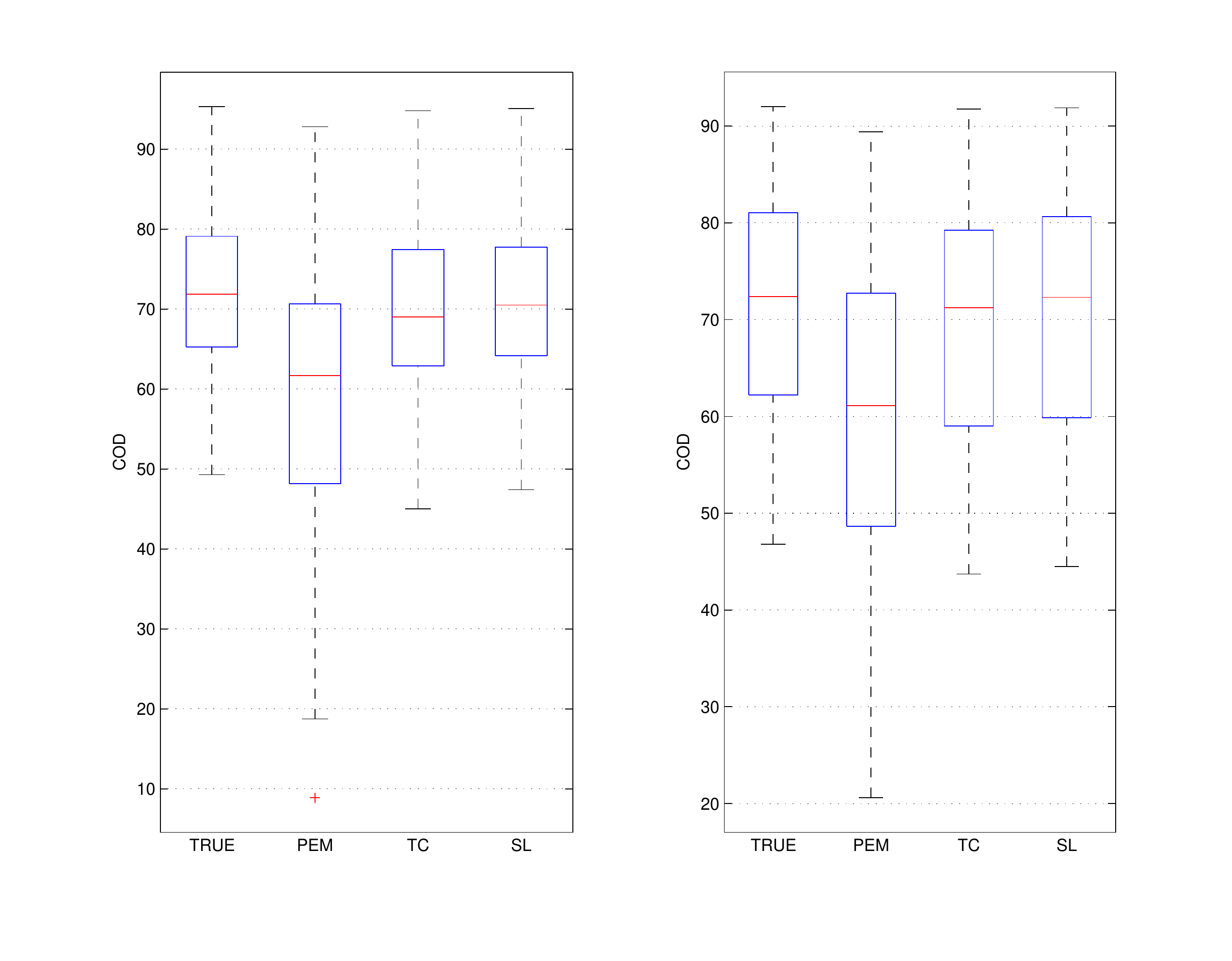}
\end{center}
\caption{One step ahead coefficient of determination in the first experiment (left panel) and in the second experiment (right panel).}\label{codSL}
\end{figure}

\begin{figure}[htbp]
 \begin{center}
\includegraphics[width=\columnwidth]{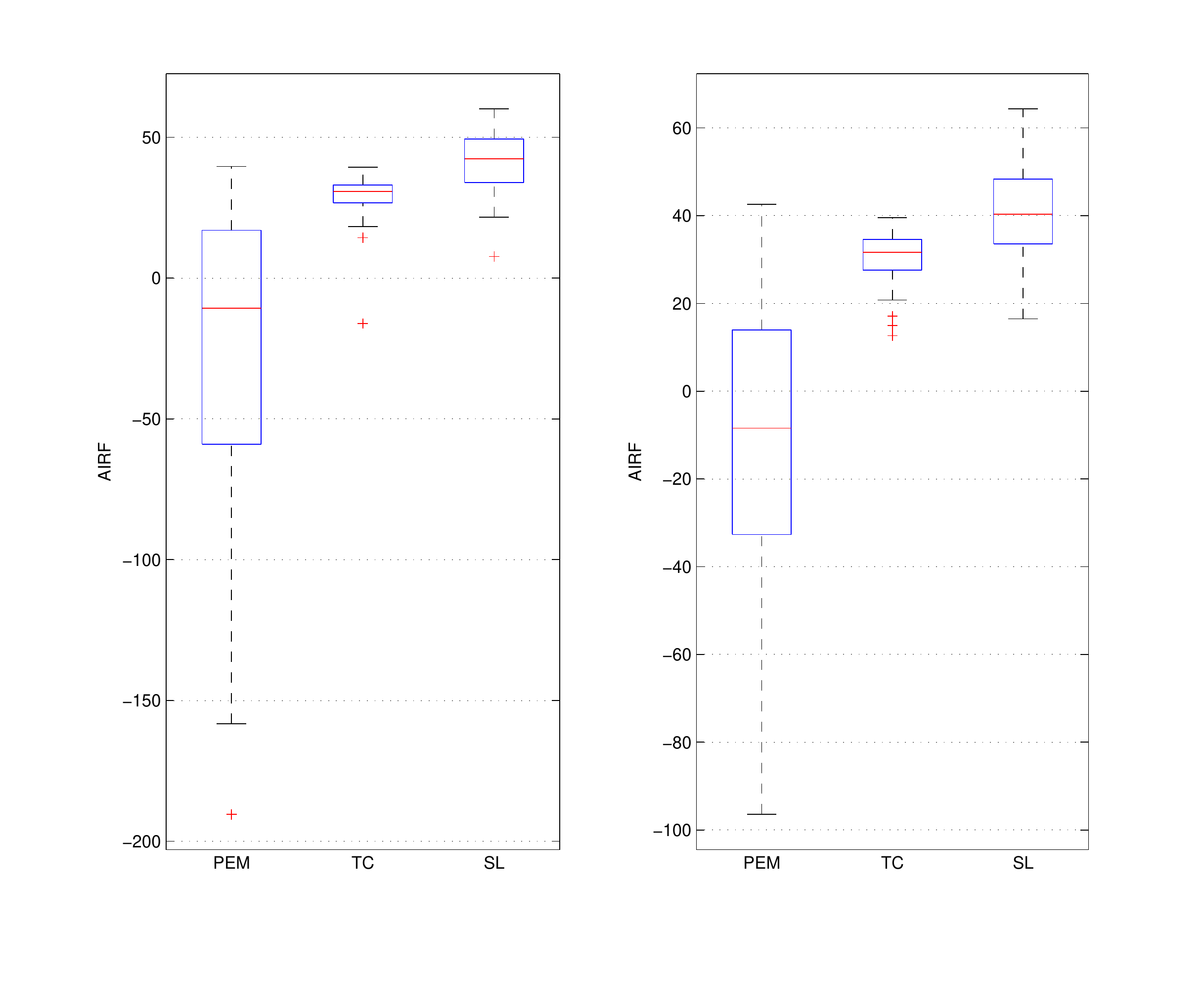}
\end{center}
\caption{Average impulse response fit in the first experiment (left panel) and in the second experiment (right panel).}\label{airfSL}
\end{figure}

\begin{figure}[htbp]
\begin{center}
\includegraphics[width=\columnwidth]{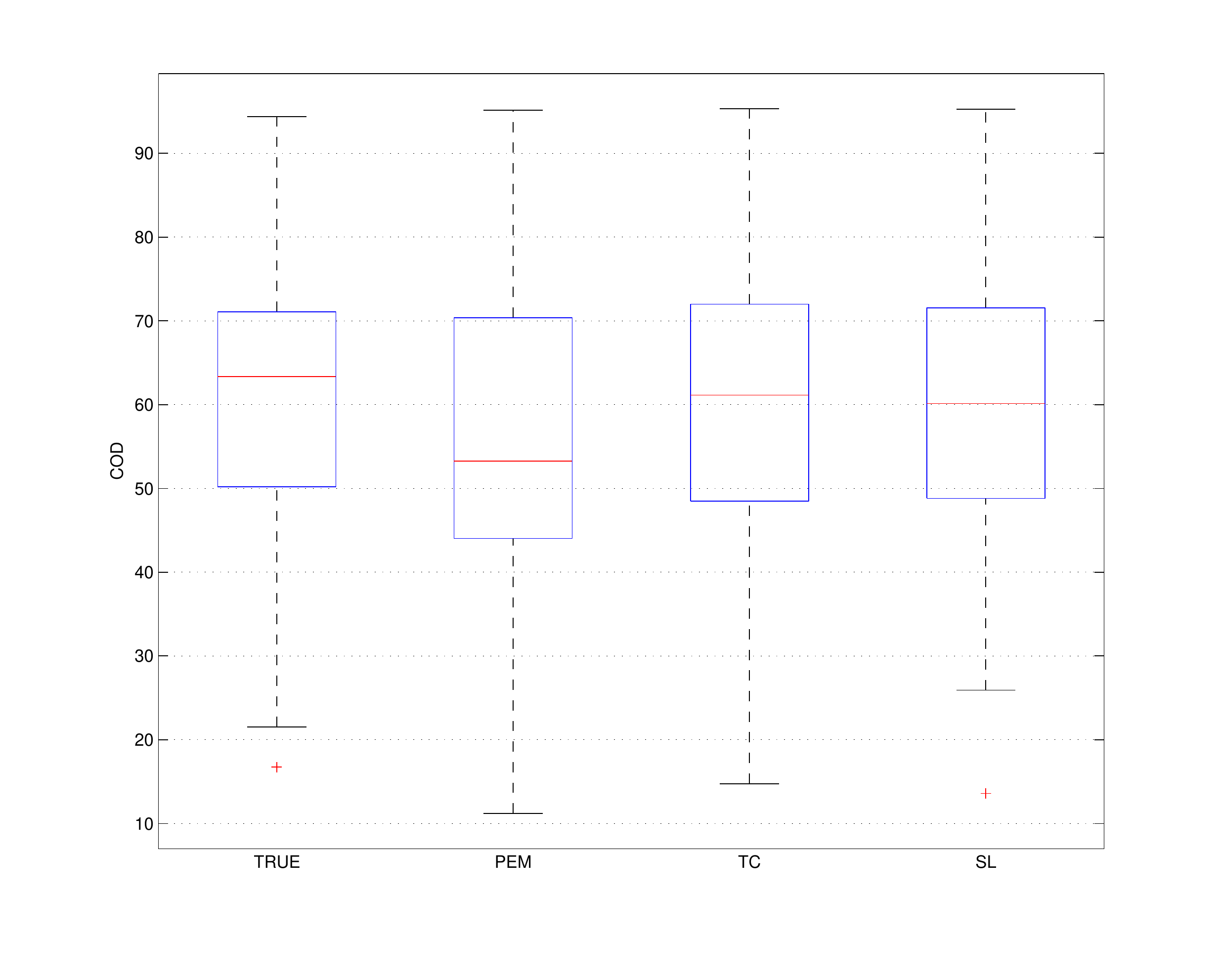}
\end{center}
\caption{One step ahead coefficient of determination in the third experiment.}\label{codgen}
\end{figure}

We compare the following one-step ahead predictors:
\begin{itemize}
  \item TRUE: this is the one computed from the true model
  \item PEM: this is the one computed from the PEM approach, as implemented in \textrm{pem.m}
function of the MATLAB System Identification Toolbox
  \item TC: this is the one computed with the approach described at the beginning of Section \ref{section_pb_formulation}
  \item SL: this is our method described in Section \ref{sec_SL_Gaussian_Regression}.
\end{itemize}
The following performance indexes are considered:
\begin{itemize}
  \item average relative complexity of the S+L network of the SL model (in percentage)
  \eq AC=\frac{100}{50}\sum_{t=1}^{50}\frac{\# SL_k}{m^2T}\eeq
    where $\# SL_k$ is the number of parameters of the estimated S+L model at the $k$-th run, whereas $m^2T$ is the number of parameters of a nonstructured model
  \item one-step ahead coefficient of determination (in percentage)
\eqn  && \mathrm{COD} =\nn\\ & & 100\left(1-\frac{\frac{1}{500}\sum_{t=1}^{500} \|y^{\mathrm{test}}(t)-\hat y^{\mathrm{test}}(t|t-1)\|^2}{\frac{1}{500}\sum_{t=1}^{500}\| y^{\mathrm{test}}(t)- \bar y^{\mathrm{test}}\|^2}\right)\nn\eeqn
where $\bar y^{\mathrm{test}}$ denotes the sample mean of the test set data $y(1)^{\mathrm{test}}\ldots y^{\mathrm{test}}(500)$ and $\hat y^{\mathrm{test}}(t| t-1)$ is the one-step ahead prediction computed using the estimated model \item average impulse response fit (in percentage)
\eq \mathrm{AIRF}=100\left(1-\frac{\sum_{k=1}^{50} \|G_k-\hat G_k\|^2}{\sum_{k=1}^{50}\| G_k- \bar G\|^2}\right)\eeq
with $\bar G=\frac{1}{50}\sum_{k=1}^{50} G_k$.
\end{itemize}
Table \ref{tabella} \begin{table}
\caption{Average relative complexity of the S+L Bayesian Network}
\label{tabella}
\begin{center}
\begin{tabular}{lccc}
\hline
Epx. \# & \#1& \#2 & \#3\\
\hline
 AC & 55.89 & 63.72 & 81.56\\
\hline
\end{tabular}
\end{center}
\end{table}  shows the percentage of the average relative complexity of the S+L network. In particular, in the first two experiments our method is able to detect that the underlying model is close to have a simple S+L network. Figure \ref{codSL} 
shows the COD in the first two experiments. One can see that SL provides a slightly better performance than TC. On the other hand, SL provides better estimators for the
predictor coefficients than the TC, Figure \ref{airfSL}. 
Finally, Figure \ref{codgen}
shows the COD in the third experiment. The median of SL is slightly worse than the one of TC. On the other hand,
the bottom whisker of SL is better than the one of TC. Indeed,  SL simplified the S+L network, see Table \ref{tabella}, increasing the robustness of the estimated
predictor impulse response coefficients.

\section{Conclusions} \label{section_conclusions}
In this paper, we proposed a Gaussian regression approach to identify multivariate stochastic processes having sparse network with few latent nodes. Simulations show
that our approach is able to identify a S+L network which does not compromise the prediction performance.

\end{document}